\documentclass[12pt,a4paper,oneside]{article}
\usepackage[top=100pt,bottom=100pt,left=100pt,right=100pt]{geometry}
\usepackage{amssymb}
\usepackage{amsfonts}
\usepackage{amsmath}
\usepackage{lastpage}
\usepackage[utf8]{inputenc}
\usepackage[bookmarks=true,colorlinks=true,linkcolor=blue,citecolor=blue]{hyperref}
\usepackage{fancyhdr}
\usepackage{color}
\usepackage{float}
\usepackage[mathlines]{lineno}
\usepackage{lscape}
\usepackage{epsfig}
\usepackage[authoryear,sort&compress]{natbib}

\setcitestyle{open={(},close={)}}
\usepackage{graphicx}
\usepackage{subcaption}
\usepackage{hyperref} 
\usepackage{caption}
\usepackage{subcaption}
\usepackage{natbib}
\bibpunct[, ]{(}{)}{;}{a}{}{,}
\usepackage{float}
\usepackage{graphicx}
\usepackage{booktabs}

\begin{document}
	\title{Modelling and Forecasting Time-Dependent Mortality Rate in Tanzania}
	\author{Samya Suleiman$^{1,2,*}$, Karl Lundengård$^{2}$, \\ John Andongwisye$^{1}$, Emmanuel Evarest$^{1}$}
	\date{}
	\maketitle
	 \noindent$^*$Corresponding author; \url{https://orcid.org/0000-0003-1713-1798} \\
	 $^{1}$ University of Dar es Salaam, Department of Mathematics, Tanzania\\$^{2}$ Imperial College of London, Department of Mathematics, London \\
	 Emails: scutesamya2804@gmail.com, k.lundengard@imperial.ac.uk, johnandongwisye@gmail.com, and sinkwembe2001@gmail.com
	\begin{quote}
	\noindent\textbf{Abstract}\\
	The increasing life expectancy enhances the importance of mortality forecasting. Most developing nations, including Tanzania, forecast mortality rates using static life tables. However, these tables exaggerate death probabilities by neglecting the declining mortality rate trend over time. This paper develops a mortality rate model based on a power exponential function, accounting for age and time variations and providing direct forecasts without methods of forecasting. The model accurately predicted mortality rates from 2023–2034 using Tanzanian data. To assess stability and reliability, a Bayesian method was employed, and numerical results revealed close alignment between the model's forecasts and Bayesian predictions.
	\end{quote}
	\begin{quote}
	 \textbf{Keywords} Mortality rates; Forecasting; Power exponential function; Bayesian method.
	\end{quote}
	\section{Introduction}
The world demographic has experienced a decrease in mortality rates, which has increased the life expectancy of individuals. Consequently, effects in planning processes for several sectors, such as insurance, demographics, health, policymaking, and pension fund systems, have occurred. Modelling and forecasting mortality rates and life expectancy are the main keys to understanding what happens in a population. Knowing this can help in solving the sector's effects by managing mortality risks and increasing the lifespan of individuals. \citep{braendle2016drives,maddison2006world}. Demographical population data with actuarial formulas can be used to create life tables, but ever since the introduction of exponential mortality models, demographic data has been smoothed and used to create improved life tables.\\

In several scientific domains, exponential mortality models are frequently employed to represent processes that display exponential growth or decay \citet{ledberg2020exponential}. Exponential models can be used to explain how the death rate evolves over time in the context of mortality rate modeling. These models allow for the possibility of a continual exponential increase or decrease in the mortality rate.\\

 Since the introduction of the Gompertz exponential mortality model \citet{gompertz1825xxiv}, which was the early mathematical model used to smooth mortality data, help to assist in extrapolation of mortality rate forecasting, create life tables, organise population information, and the development of the countries, But it had the weakness of overestimating death rates at ages older than 80.\\

Some researchers developed mortality rate models by extending the Gompertz model. \citet{makeham1860law}, \citet{thiele1871mathematical},  \citet{weibull1951statistical} and \citet{perks1932some} developed mortality rate models that provided a satisfactory fit to adult mortality but overestimated death rates at ages below 20 and older than 80. Recently, models by \citet{hannerz2006}, \citet{skiadas2018fokker}, and \citet{lundengard2019}, came up with models that give a comparatively good fit to mortality rates over the entire age range by considering the number of parameters and model complexity to follow the complete pattern of human life.\\
    
      Improvements in good healthcare and hygiene, resulting from the introduction of good medicine capable of reducing death, have disturbed the economic security of nations and the density and structure of their populations. This has forced researchers to do forecasting of mortality rates due to the increase in people's life expectancy \citep{pascariu2018modelling}. \citet{bennett2015future} demonstrated that a better statistical and stochastic model for mortality rate forecasting can encourage a country to develop a good quality of life as it will improve the economy of the respective country. Over the last few decades, many methods have been developed for forecasting mortality using stochastic models, such as \citet{alho1992modeling}, \citet{bell1997comparing}, and \citet{lee1992modeling}.\\

	  \cite{Belinda2018comparison}, \citet{sulemana2019comparison}, \citet{lundengard2019}, and \citet{boulougari2019application} used power exponential models that vary with only ages in forecasting mortality rates using the Lee and Carter Method (LCM). \cite{bongaarts2005long} established and put to the test a new version of the logistic model, which described how the pattern of mortality rates changes in both ages and time. But then we directly predicted future trends only in adult mortality rates for ages above 25 and compared the results by forecasting with the LCM.\\

	   Even though the method has been applied in different countries such as Sri Lanka, India, and Thailand \citep{aberathna2014modeling, chavhan2016modeling, yasungnoen2016forecasting} and others, it has some weaknesses. It ignores the error of the estimated parameters and accounts for the error of the mortality index by assuming the forecast mortality index error dominates all others. Leads to prediction intervals that are too narrow. Furthermore, \citet{lee1992modeling} stated that, for the method to give accurate results of mortality rate forecasting, it needs a long range of past mortality rate data. These were the reasons for extending the LCM to state space models (SSM).\\

	          In SSM, all parameters are jointly estimated, and their uncertainties can be calculated. \citet{pedroza2002bayesian}, \citet{fung2015state}, \citet{pedroza2006bayesian}, and \citet{kogure2010bayesian} applied SSM with the Bayesian method to forecast mortality rates. This method has been used for jointly modelling time series and incorporating all-parameter uncertainty. \citet{booth2008mortality} said the principal of SSM is that the statistical and characteristic procedures of the entire estimation and forecasting process are well familiarised and direct. This method has more advantages than LCM, which employs stages in parameter estimation and fails to integrate uncertainties.\\

	  Despite several studies done on models related to mortality rates, models that vary with both ages and time are very rare. Based on the model by \citet{lundengard2019} that varies with only ages, this study has developed a new mortality rate model that describes how mortality rates vary with both ages and time. Furthermore, direct forecast mortality rates for a whole age range with limited Tanzanian data. The results were compared with the Bayesian method.\\

	     This paper is organised as follows: Section \ref{sec:material_and_methods} gives an overview of the data, mathematical model, and methods used in this paper. Findings and results are discussed in Section \ref{sec:results}. In Section \ref{sec:conclusion}, we give conclusions and remarks and discuss future work.      
	
	\section{Materials and Methods}
	\label{sec:material_and_methods}
	\subsection{Materials}
Tanzania regional-wise population projection data for 2012, \citet{nbs2002census} and population projections for 2013–2035, \citet{nbs2013projection}  from NBS reports with single-ages from 0–78+ were used to generate mortality rates for the period of 2012–2034. \\     

\begin{figure}[H]
		\centering
		\includegraphics[width=0.8\textwidth]{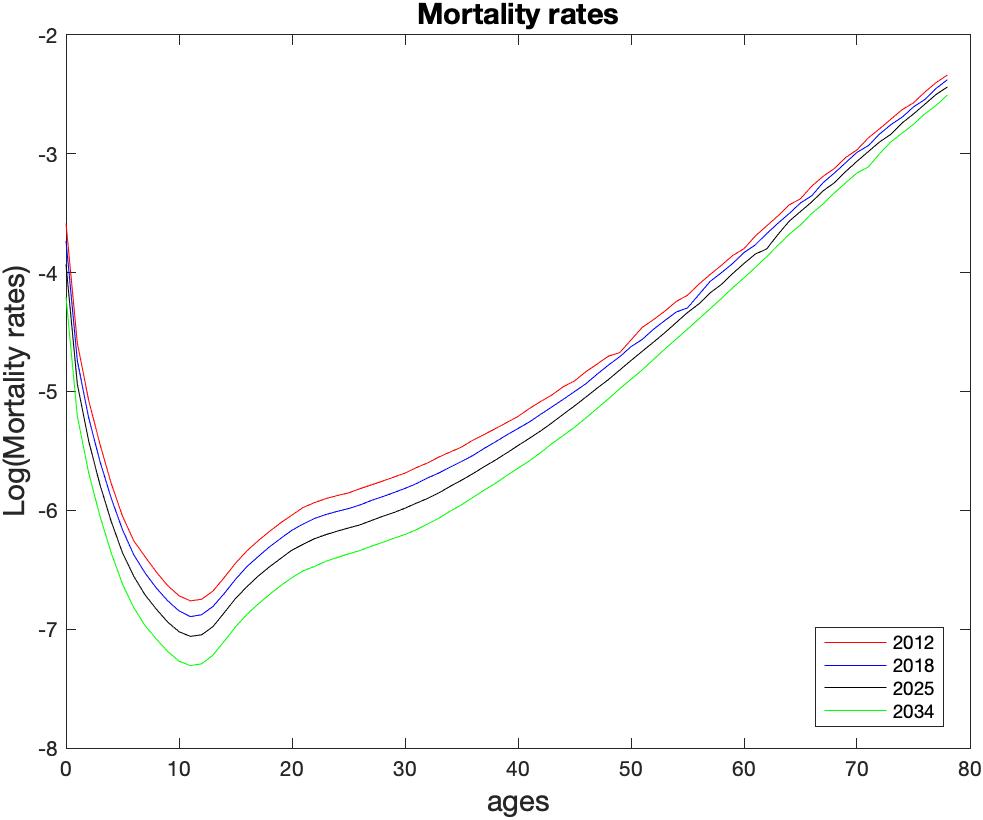}
		\caption{Generated Tanzania mortality rates for Years 2012, 2018, 2025 and 2034.} \label{fig:mortality} 
	\end{figure}

The mortality rates generated were calculated using the formula given by:
	\begin{equation}\label{eq:mortalityrates}
		\mu_{x,t} = \frac{d_{x,t}}{l_{x,t}},
	\end{equation}
    where, $\mu_{x,t} $ is the mortality rate of individuals aged $x$ in year $t$. $d_{x,t}=l_{x,t}-l_{x,t+1}$ is the number of deaths of individuals aged $x$ in year $t$, and $l_{x,t}$ is the number of people alive at age $x$ in year $t$.\\ Figure~\ref{fig:mortality}  represents  the  logarithm  of  generated  mortality  rates data for Tanzania for a few  selected years 2012, 2018, 2025, and 2034.\\
         
    The shapes in the graph show a decrease in mortality rates over the years that is caused by improvements in good healthcare and hygiene in the country. The shapes in Figure~\ref{fig:mortality} also indicate that, mortality rates follow the same pattern and trend as mortality rates in other countries; see other examples in \citep{boulougari2019application}.


	\subsection{Methods}
	\subsubsection{Mortality rate model formulation.}
	Consider the mortality rate  model $\mu_{x}$ by  \citet{lundengard2019} as follows, 
	\begin{equation} \label{eq:model}
		\mu_{x}= \frac {a_{1}e^{a_{2}x}}{x} + b_{1} \left(x e^{-b_{2}x}\right)^{b_{3}}.
	\end{equation}
    Parameters, $a_{1}$ and $a_{2}$ represent mortality rates for ages before adolescence and old age, respectively. While parameters $ b_{1},b _{2}$ and $b_{3}$ describe the effect of the hump in the mortality rates for ages between adolescence and early adulthood around 15–25 years, as shown in Figure~\ref{fig:mortality}. This study extends the model from Equation~\ref{eq:model} to include both ages $x$ and year $t$. The mortality rate that varies with both ages and time is defined as follows:
	
		\begin{equation} \label{eq:extended}
		\hat{\mu}_{x,t}= \frac{a_1(t)e^{a_2(t)x}}{x } + b_1(t) \left(x e^{-b_2 (t)x}\right)^{b_3(t)}.
	\end{equation}
     \noindent where, parameters are defined in similar way as in Equation~\ref{eq:model} with the addition to change with time trend.\\
     
	The formulation of the model of Equation \ref{eq:extended} was done by extending the parameters of Equation~\ref{eq:model} to vary with time. To manage that, the generated mortality rates were fitted to model Equation~\ref{eq:model}, model parameters were estimated. Then, the analysis of the estimated parameters was done by studying the behaviour of how the parameter values changed over time.

	\subsubsection{Model Parameter estimations.}
Thereafter, the values of $a_{1}(t), a_{2}(t), b_{1}(t), b_{2}(t)$ and $b_{3}(t)$ were estimated using the non-least squares method by minimising the sum of squared residuals $R$ given by the following formula:

\begin{equation} \label{eq:NLLM}
R=\min\left(\sum_{x=0}^{n}{r_{x}}^{2}(t)\right) ,
\end{equation}

 \noindent where  $n$ is the maximum age, $r_{x}(t)= \mu_{x,t}-\hat{\mu}_{x,t}$ with $ \mu_{x,t}$ and $\hat{\mu}_{x,t}$ be mortality rates from Equation \ref{eq:mortalityrates}  and Equation \ref{eq:extended} respectively. The equation \ref{eq:mortalityrates}  has a set of data points of $x=0,2,\cdots,78+$ and years $t=2012,\ldots,2034$. The parameters were estimated from year $2012$ and were estimated with changing initial parameter values, that is, the obtained parameter of year $t$ was used to estimate the parameter of the next year $t+1$, and so on.

	\subsubsection{Forecasting with Bayesian Method}
	 To check the reliability and stability of the model, the Bayesian method was employed based on  \citep{fung2015state} . It is a statistical method that offers results concerning unknown parameters of SSM using probability statements \citep{pedroza2002bayesian}. Parameters are summarised by an entire distribution of values rather than one mounted value as in classical frequent analysis as has been in \citep{lee1992modeling}.\\
	
	 The estimation of parameters was done using the posterior distribution of each parameter, depending on the chosen prior and likelihood of the parameters. Then, the posterior distributions of each parameter are accessible by sampling using the Kalman filter method and Markov Chain Monte Carlo (MCMC) methods. Consider the SSM defined in Equation~\ref{eq:statespacemodel}
	
	\begin{equation} \label{eq:statespacemodel}
		y_{t} = \alpha_{x} + \beta_{x} k_{t} + \varepsilon_{t}  \qquad  \textrm{such that}  \qquad k_{t} = k_{t-1} + \lambda + \omega_{t}.
	\end{equation}

	 \noindent $y_{t}$ is the generated mortality rate in year $t$ as $\mu_{x,t}$. $\alpha_{x}$ is interpreted as the average pattern of mortality at age $x$ which describes the general shape of the mortality curve. $\beta_{x} $ represents how mortality at each age varies when the general level of mortality changes that is, it represents each age group's response or susceptibility to the temporal mortality index. $k_t$  is the temporal mortality index that captures the evolution of rates over time. $\lambda$  is a drift parameter that measures the constant annual change in the series of $k_t.$  $\varepsilon_{t} \overset{iid}\sim  N(0,\sigma_{\varepsilon}^{2}1i) \hspace{0.1cm} \text{and} \hspace{0.1cm} \omega_{t} \sim N(0, \sigma_{\omega}^{2})$ represent error terms that cause the deviation of the model.  $1i$ is the $i$ by $i$ identity matrix which makes $\varepsilon_{t}$ a diagonal matrix, and $N(\cdots)$ is a Gaussian distribution, then for $i = 1\cdots x$. \\

	The law of conditional probability is applied to provide probabilistic inferences concerning any parameter \citep{west1997}. Let $\phi$ be the set of all parameters $\alpha_{x}$,$\beta_{x}$, $\lambda$, $\sigma_{\varepsilon}^{2}$ and $\sigma_{\omega}^{2}$ and $k_{t}$ is a state parameter.
	The conditional probability  $p({\phi}$/${y_{t}})$ is defined by; 
	$$p({\phi}/{y_{t}})= \frac{p({y_{t}}/{\phi}) \cdot p(\phi)}{p(y_{t})}$$ 
	where  $p(y_{t})= \int_{\phi}{p({y_{t}}/{\phi})\cdotp(\phi)}.$ Since $y_{t} =\mu(x,t)$ are generated mortality rate values then $p(y_{t})$ is constant. The proportional form of Bayesian theorem becomes, 
	\begin{align*} 		
		\begin{array}{ccccc}
			p({\phi}/{y_{t}} )\propto & & p({y_{t}}/{\phi})  \cdot  & & p({\phi}) \\       
			\downarrow & & \downarrow & & \downarrow \\
			\text{posterior density}  & & \text{likelihood} & &\text{prior density}
		\end{array} 
	\end{align*}
	Then posterior density $p(k_{t},\phi /y_{t})$ of the state $k_{t}$ as well as parameters $\phi$ given the observations $y_{t}$ are obtained using Kalman filter and MCMC respectively.\\
	
    The Kalman Filter is an optimal estimator consisting of a set of mathematical equations that provides a recursive computational methodology for estimating the state of a discrete--data. The estimation of state $k_{t}$ is recursively performed using two Kalman filter steps \citep{petris2009dynamic}. \\First, we find the mean $m_{t}$ and variance $C_{t}$ of $N(m_{t},C_{t})$ by assuming the initial values of $m_{t-1}$ and $C_{t-1}$. Then find the distribution of parameter values of state $k_{t}$ which are distributed normally with mean and variance given by $ N(h_{t},H_{t})$ in step two as described as follows: \\

	Step one comprises the  update steps (a), (b) and (c). Then (d) predict state parameter $k_{n}$ as follows:\\
	(a) Posterior for $k_{t-1}$: $(k_{t-1} | y_{t-1}) \sim N(m_{t-1}, C_{t-1})$.\\
	(b) Prior for $k_{t}$: $(k_{t} | y_{t-1}) \sim N(a_{t},R_{t}), \text{where} \hspace{0.1cm}a_{t}= m_{t-1}+ \lambda , R_{t}= C_{t-1}+ \sigma_{\omega}^{2}$.\\
	(c) 1-step forecast: $(y_{t} | y_{t-1}) \sim N(f_{t},Q_{t}), \text{where} \hspace{0.1cm}f_{t}=\alpha+\beta a_{t} , Q_{t}=\beta R_{t}\beta^{'}+ \sigma_{\varepsilon}^{2}1_{i}$.\\
	(d) Posterior for $k_{n}$: $(k_{t} | y_{t}) \sim N(m_{t},C_{t}), \text{where} \hspace{0.1cm} m_{t}= a_{t}+ s_{t}(y_{t}- f_{t}),\hspace*{0.1cm} C_{t}= R_{t}(1- s_{t}\beta) , \text{where} \hspace*{0.2cm}s_{t}= R_{t}\beta^{'}Q_{t}^{-1}.$ 
	
	Step two;\\
	After obtaining $k_{n}$ then we do backward sampling from $t = n-1, n-2, . . . , 1, 0 $ to get $k_{t}$ with $ N(h_{t},H_{t})$ as $h_{t}= m_{t}+C_{t}R^{-1}_{t+1}(k_{t+1}- a_{t+1}) \hspace{0.1cm} \text{and} \hspace{0.1cm} H_{t}= C_{t} - C_{t}R^{-1}_{t+1}C_{t}$. \\ Where $R_{t+1}=C_{t}+\sigma_{\omega}^{2}$. Then $k_{t}$ is approximated as: 
	\begin{equation}\label{eq:9}
		k_{t} \sim \left( h_{t}= m_{t}+C_{t}R^{-1}_{t+1}(k_{t+1}- a_{t+1}), H_{t}= C_{t} - C_{t}R^{-1}_{t+1}C_{t}\right)
	\end{equation}
	
	To estimate parameters $\phi$, the conditional distributions of each parameter is obtained from posterior distribution $p(\phi/k_{t}/y_{t})$ and MCMC sampling is used to sample the parameters.\\
	 Vague prior values were used as in \citet{fung2015state} for parameters $\alpha_{x}$, $\beta_{x}$, $\lambda$, $\sigma_{\varepsilon}^{2}$ and $\sigma_{\omega}^{2}$ which are $\alpha_{x} \sim N(\hat{\mu}_{\alpha} ,  \hat{\sigma}_{\alpha}^{2} )$ , $\beta_{x} \sim N(\hat{\mu}_{\beta} ,  \hat{\sigma}_{\beta}^{2})$ , $\lambda  \sim N(\hat{\mu}_{\lambda},  \hat{\sigma}_{\lambda}^{2})$, $\sigma_{\varepsilon}^{2} \sim IG(\hat{a}_{\varepsilon},  \hat{b}_{\varepsilon})$ and $\sigma_{\omega}^{2} \sim IG(\hat{a}_{\omega},  \hat{b}_{\omega})$ where IG means Inverse Gamma distribution. It is assumed that the prior for all parameters are independent. Then, posterior distribution of the parameters shares the same distribution type as the prior distributions.\\
	\\
	 MCMC simulations allow for parameter estimation such as means, variances, expected values, and exploration of the posterior distribution of Bayesian models \citep{kavetski2006bayesian}. To assess the properties of posterior distributions, many representative random values should be sampled from the distribution using MCMC \citep{mbalawata2014adaptive}.\\       
	    \\
	     This study used Gibbs sampling to find final estimations for the posterior distributions since it chose a random value for a single parameter, holding all the other parameters constant. Then samples obtained from MCMC Gibbs sampling immediately represent samples from the posterior distribution for forecasting. The sampling was done using MATLAB for 12000 samples and burn-in 8000. The last 4000 samples were selected to get parameters of $\phi$ that gave reasonable results for each parameter. \\       
	     \\
	      The length of the burn--in period was judged visually by plotting sample values from the MCMC process. Assuming the total number of samples remained after burn--in is $j$, and $Z_{\phi}$ is the set of samples for each parameter,  for  $g= 1,2,\ldots,j$ the mean of samples for all parameters was calculated as:

	\begin{equation} \label{eq:meanvalue}
		\text{Mean of the sample}(Z_{\phi}) = \frac{1}{g}{\sum_{i=1}^{g}{\phi}^{i}}
	\end{equation}  
	
	\noindent The mortality index  used to forecast mortality rate with SSM was given by;
	\begin{equation} \label{eq:mortalityindex}
		k_{t+h}^{(g)} \sim N\left(k_{t+h-1}^{(g)}+\lambda^{(g)}, (\sigma_{\omega}^{2})^{(g)}\right), 
	\end{equation}
	
	\noindent where $h$, is the forecast horizon with $h \geq 1$. Finally, the mortality index from Equation~\ref{eq:mortalityindex} with estimated parameters from MCMC were used to obtain the forecasted mortality rates using the following equation,  
	\begin{equation} \label{eq:forecastmortality}
		y_{t+h}^{(g)} \sim N\left(\alpha^{(g)}+ \beta^{(g)}k_{t+h}^{(g)}, (\sigma_{\epsilon}^{2})^{(g)}I_{p}\right).
	\end{equation}

	\newpage
	\section{Results and Discussion.} \label{sec:results}
\begin{figure}[H]
	\centering
	\includegraphics[width=0.8\textwidth]{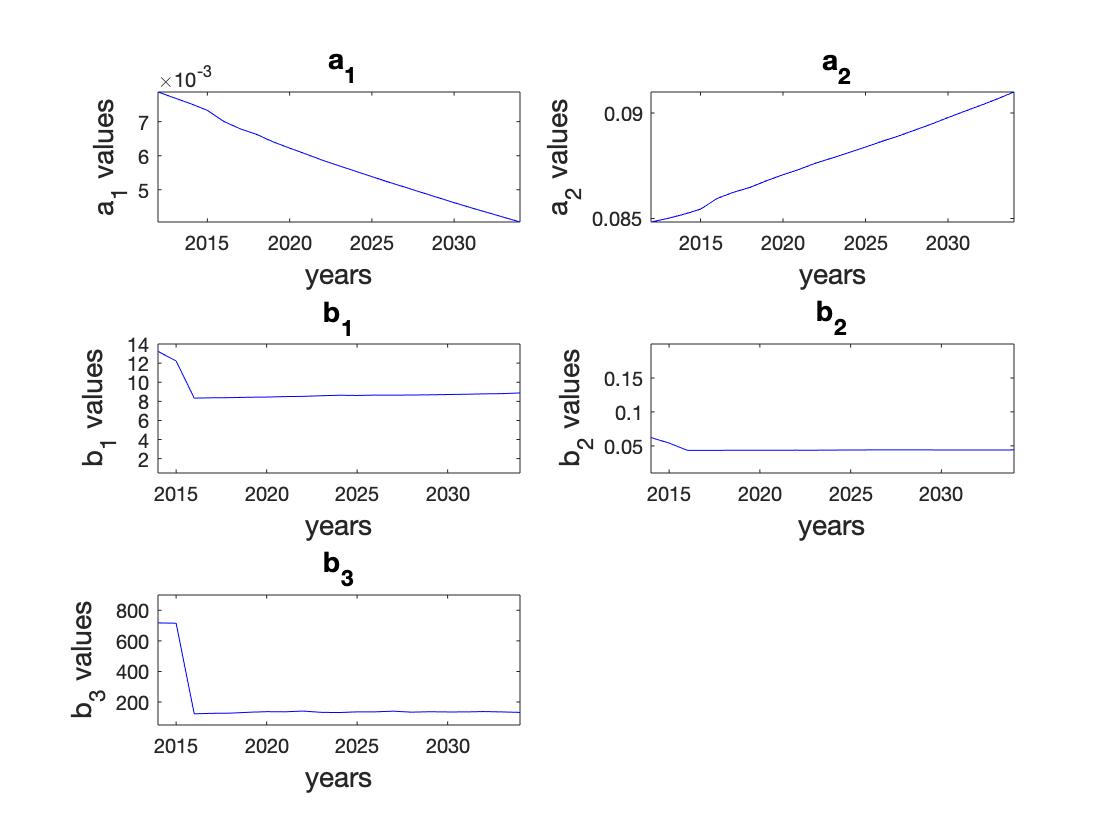}
	\caption{Estimated parameters in relation to time in years.}
	\label{fig:parameters} 
\end{figure}

 The study first  fitted generated mortality rates to the model from Equation~\ref{eq:model} to estimate the parameters. Then, the analysis of parameter values was done by studying the behaviour of how the estimated parameters change over time. The shapes of the parameters of the model appeared to change exponentially for parameters $a_{1}$ and $a_{2}$ and change linearly with $b_{1},$ $b_{2},$ and $b_{3}$ as shown in Figure~ \ref{fig:parameters} \\
 
The shapes and trend of parameters in relation to time, shown in  Figure \ref{fig:parameters} led the modification of parameters of the Equations \ref{eq:model} to be as shown in the following equations 
\begin{align}
	\label{eq:a_1} 
	a_{1}(t) & = a_{1}^{0}e^{-K_{1}t}, \\
	\label{eq:a_2}
	a_{2}(t) & = a_{2}^{0}e^{K_{2}t}, \\
	\label{eq:b_1}
	b_{1}(t) & = b_{1}^{0} + K_{3}t, \\
	\label{eq:b_2}
	b_{2}(t) & = b_{2}^{0} + K_{4}t, \\
	\label{eq:b_3}
	b_{3}(t) & = b_{3}^{0} + K_{5}t.
\end{align}

The parameters $a_{1}^{0}$, $a_{2}^0$, $b_{1}^{0}$, $b_{2}^{0}$, $b_{3}^{0}$ are chosen values that give initial  values of Equations~\ref{eq:a_1}--\ref{eq:b_3} at  $t=0$. The $K_{1}$, $K_{2}$, $K_{3}$, $K_{4}$ and $K_{5}$  are the slopes of the parameters. The values of $K_{1}$ and $K_{2}$ are obtained by finding the first derivative of $a_{1}(t)$ and $a_{2}(t)$ with respect to $t$, while $K_{3}$, $K_{4}$ and $K_{5}$ are obtained using the general equation for finding the slope of a line. The value of $a_{1}$ tends to decrease as mortality rates at low ages decrease. 

To estimate parameter values  $a_{1}(t), a_{2}(t), b_{1}(t), b_{2}(t)$ and $b_{3}(t)$ in Equation~\ref{eq:extended}, we first estimated parameters  of Equation~\ref{eq:model} using the least square method with generated data. Values of $a_{1}^{0}$, $a_{2}^0$, $b_{1}^{0}$, $b_{2}^{0}$, $b_{3}^{0}$, $K_{1}$, $K_{2}$, $K_{3}$, $K_{4}$ and $K_{5}$ were chosen as shown in Table~\ref{table:1}  so that they match well the estimated values for $a_{1}(t), a_{2}(t), b_{1}(t), b_{2}(t)$ and $b_{3}(t)$. \\

	\begin{table}[H]
	\centering
	\begin{tabular}{llcll}
		Parameter   & Value      & & Parameter & Value                  \\ \cmidrule{1-2} \cmidrule{4-5}
		$a^{0}_{1}$ & $0.006223$ & & $K_{1}$   & $3.4631 \times10^{-2}$ \\ \cmidrule{1-2} \cmidrule{4-5}
		$a^{0}_{2}$ & $0.08707$  & & $K_{2}$   & $3.1204 \times10^{-3}$ \\ \cmidrule{1-2} \cmidrule{4-5}
		$b^{0}_{1}$ & $8.446$    & & $K_{3}$   & $3.0111 \times10^{-2}$ \\ \cmidrule{1-2} \cmidrule{4-5}
		$b^{0}_{2}$ & $0.0434$   & & $K_{4}$   & $3.5 \times10^{-5}$    \\ \cmidrule{1-2} \cmidrule{4-5}
		$b^{0}_{3}$ & $132.5$    & & $K_{5}$   & $5.22 \times10^{-1}$   \\ \cmidrule{1-2} \cmidrule{4-5}
	\end{tabular}
	\caption{Shows the values of parameters chosen to give the initial values in Equations \ref{eq:a_1}--\ref{eq:b_3} for the developed model when $t=0$.}
	\label{table:1} 
\end{table}

\begin{figure}[H]
	\centering
	\includegraphics[width=0.8\textwidth]{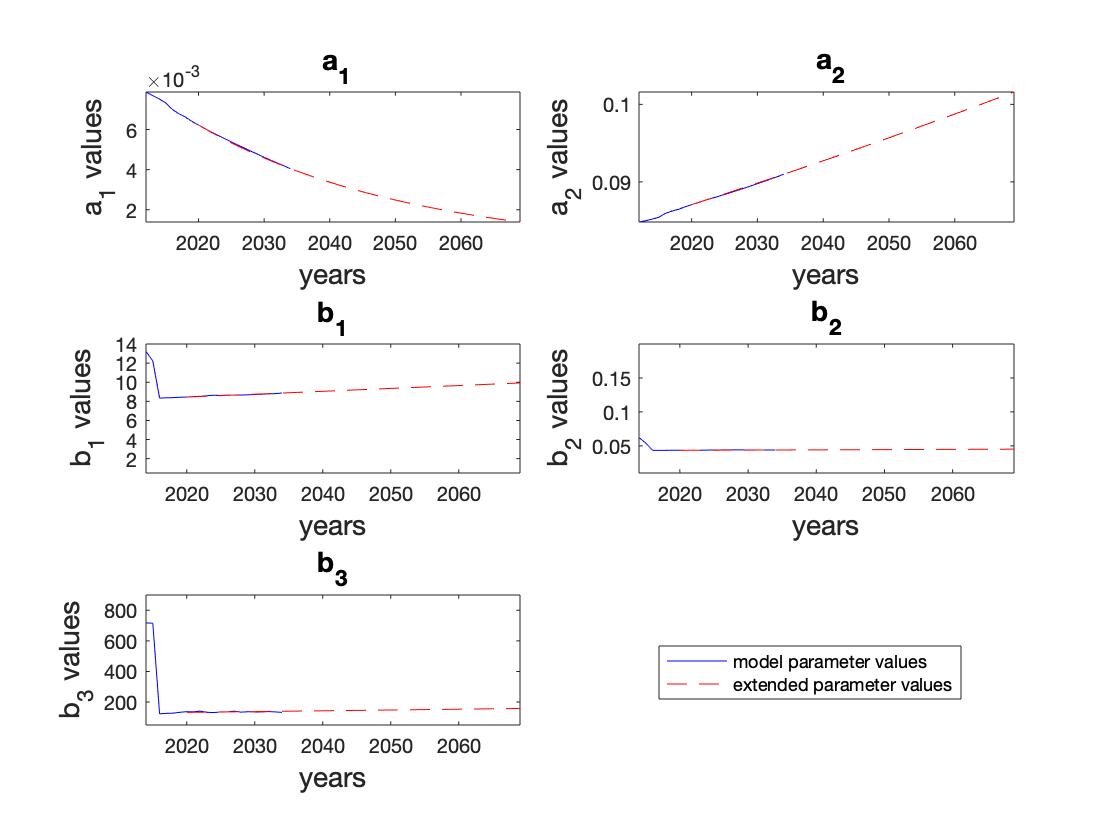}
	\caption{The overlapping of parameter values with time in years.}
	\label{fig:overlapping} 
\end{figure}

 Using Equation~\ref{eq:NLLM} we estimated values of parameters $a_{1}(t), a_{2}(t), b_{1}(t), b_{2}(t)$ and $b_{3}(t)$  with changing  initial values of parameter for each year from 2012--2070, that is, the obtained parameter of year 2012 was used to estimate parameter of  2013, and so on. Then plotted from 2020--2070 with parameters of Equation~\ref{eq:model}. The parameters appeared to overlap each other as shown in Figure \ref{fig:overlapping}. The overlapping behavoiur from 2020--2034 made the study to validate and enhance the extension of model parameters 
 to Equation~\ref{eq:extended}.\\
 
\begin{figure}[H]
	\centering
	\includegraphics[width=0.8\textwidth]{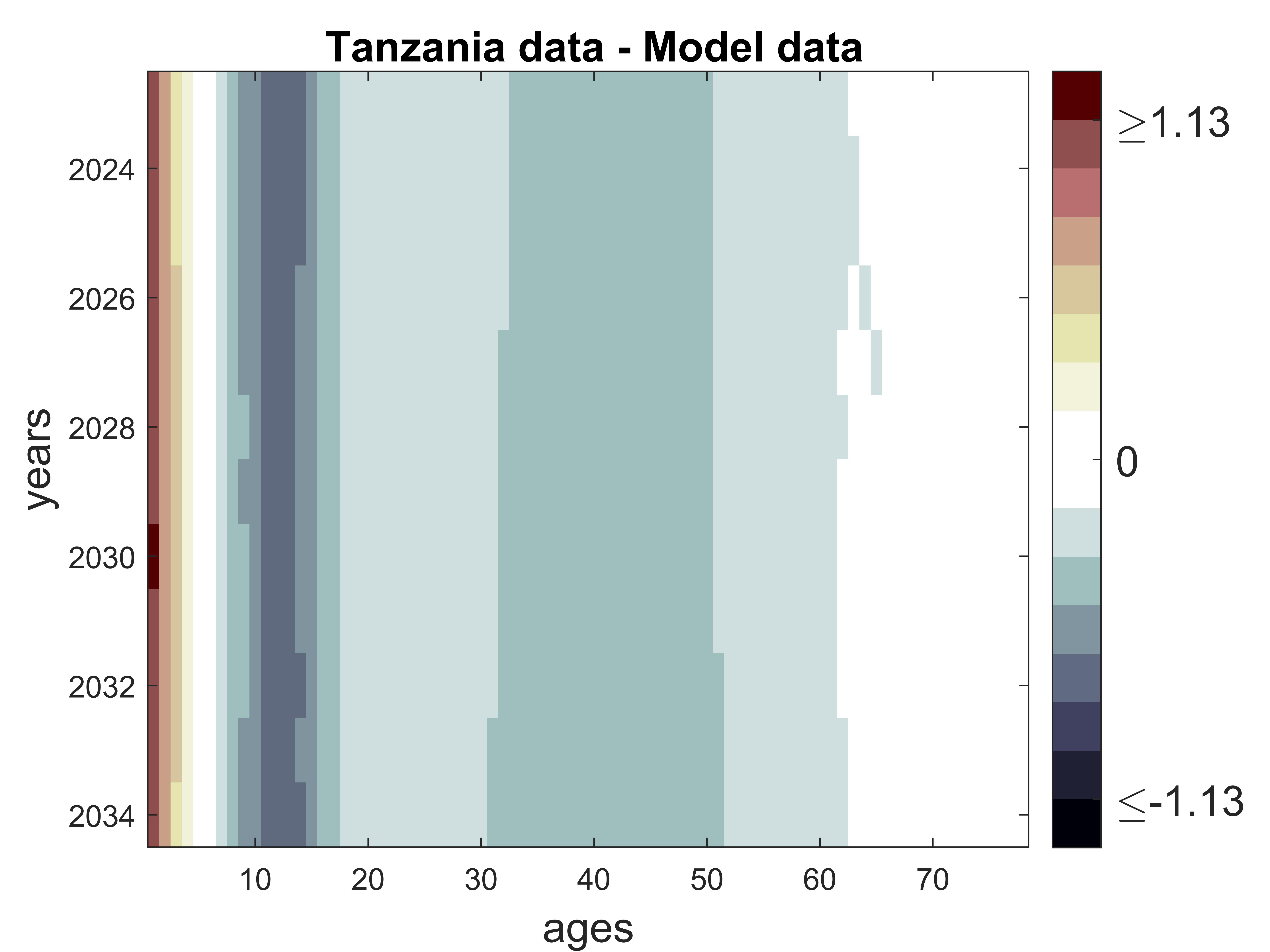}
	\caption{The difference in mortality between Tanzania data  and forecast results from the developed model for the year 2023--2034.}
	\label{fig:forecastmodel_difference} 
	\end{figure}

    \noindent Then, a new set of data was obtained by changing the initial parameters for the years 2012–2070. By using the extended model, mortality rates were then forecasted directly without using the method of forecasting for 12 years from 2023–2034. Model validation was done by comparing the forecast results from the developed model with generated data for years 2023--2034 as shown in Figure \ref{fig:forecastmodel_difference}. The mortality rate forecasting was done for 2034 due to the available data. \\
\\
	
	Also, in Bayesian forecasting, Tanzania generated data for the years 2012--2022 as past data was fitted to the method. Then assumed initial parameters by \citet{fung2015state} for $\alpha_{x}$, $\beta_{x}$, $\lambda$, $\sigma_{\varepsilon}^{2}$ and $\sigma_{\omega}^{2}$ to be $-5, 0.2, -0.8,$ $2\times 10^{-1}$ and $6.5\times 10^{-2}$ respectively from Equation \ref{eq:statespacemodel} were used to generate samples and estimate posterior distributions for both $k_{t}$ and parameters $\phi$. \\
	
	To estimate parameters $\alpha_{x}$, $\beta_{x}$, $\lambda$, $\sigma_{\varepsilon}^{2}$ and $\sigma_{\omega}^{2}$  of SSM with Bayesian method, we used initial priors of state $k_{t}$ in step one of Kalman Filter of $m_{0}$ and $C_{0}$ to be 0 and 100 respectively. Again initial priors of $\phi$ as $\mu_{\alpha_{x}},$ $\mu_{\beta_{x}},$ and $\mu_{\lambda},$ to be 0, also $\sigma_{\alpha_{x}}^{2},$ $\sigma_{\beta_{x}}^{2},$ and $\sigma_{\lambda_{x}}^{2} $ to be 100 while ${a}_{\epsilon},$ ${a}_{\omega}$ to be 2.1 and ${b}_{\epsilon},$ $ {b}_{\omega}$ to be 0.3. Then mortality rates were forecasted for a period of 2023--2034. \\
	
	Figure~\ref{fig:forecastmodel_difference} below shows the difference in future mortality rates between generated Tanzania data and forecasts from developed model.\\
	
 The values on the right side of the figure show the increasing and decreasing mortality rate. Negative values (blue) indicate that the forecast from the developed model produces a higher mortality rate than that of the generated data. Apart from that, positive values (red) indicate that the forecast from the developed model produces a lower mortality rate than that of the generated data. Also, the white  colour indicated that the forecasts are very similar. Generally, it can be noted that the greater the difference in forecast, the darker the corresponding colour.\\
 \\
The forecast from the developed model appears to have low mortality rates at lower ages particularly below the age of 5 years, and high mortality rates at ages near the hump compared to the generated data. This suggests that the model can be underestimating and overestimating mortality at lower ages before 15. Also, it shows close matches to the generated data for the ages between 18–30 as the colour is almost going to white. Furthermore, it also shows reasonably accurate results for individuals aged 65 and above.\\
\\

\begin{figure}[H]
\centering
\includegraphics[width=0.8\textwidth]{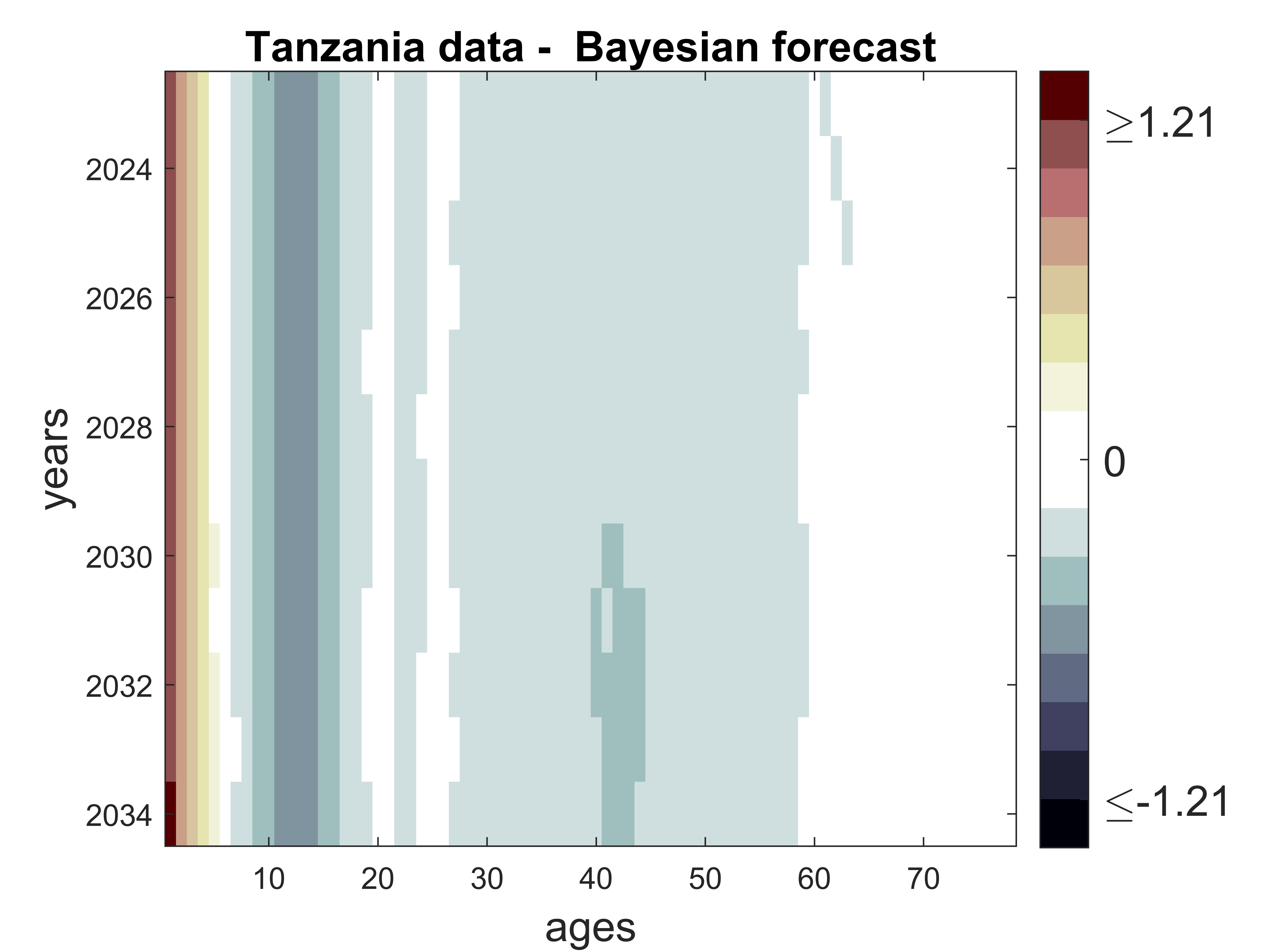}
\caption{ The difference between Tanzania generated data and forecast results from Bayesian method for the year 2023--2034. }
\label{fig:forecastBayesian_difference}
\end{figure}

The difference in forecast between the Bayesian method and generated data is shown in the Figure~\ref{fig:forecastBayesian_difference}, gives almost similar behavior when forecasting with the developed model.\\

The Bayesian method and model developed showed a similar forecast in mortality rates for ages between 0–15 with Tanzanian data. However, the Bayesian method performed better with Tanzanian data around 20–25 years of age compared with the forecast of the developed model. But both are more effective at ages above 60. This suggests the utility of the Bayesian approach and the developed model in providing more accurate mortality rates at higher ages. Around 30–50 Bayesian forecast values are lower compared to the developed model when forecasted with Tanzanian data. This indicates the effectiveness of the Bayesian method during those ages.

\begin{figure}[H]
\centering
\includegraphics[width=0.8\textwidth]{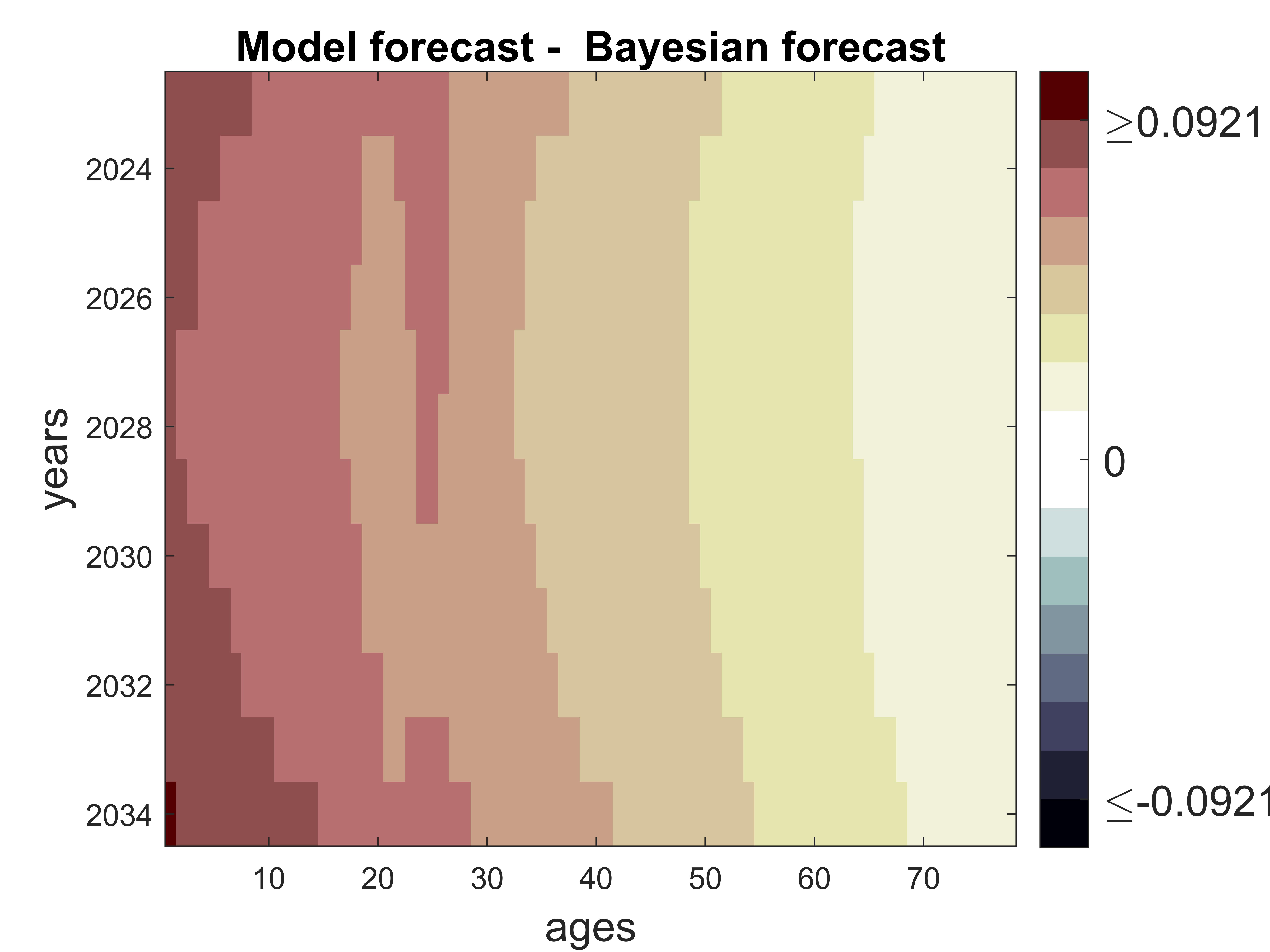}
\caption{The difference between forecast results from the developed model and Bayesian method.} 
\label{fig:forecastmethods_difference}
\end{figure}

\begin{figure}[H]
\centering
\includegraphics[width=0.8\textwidth]{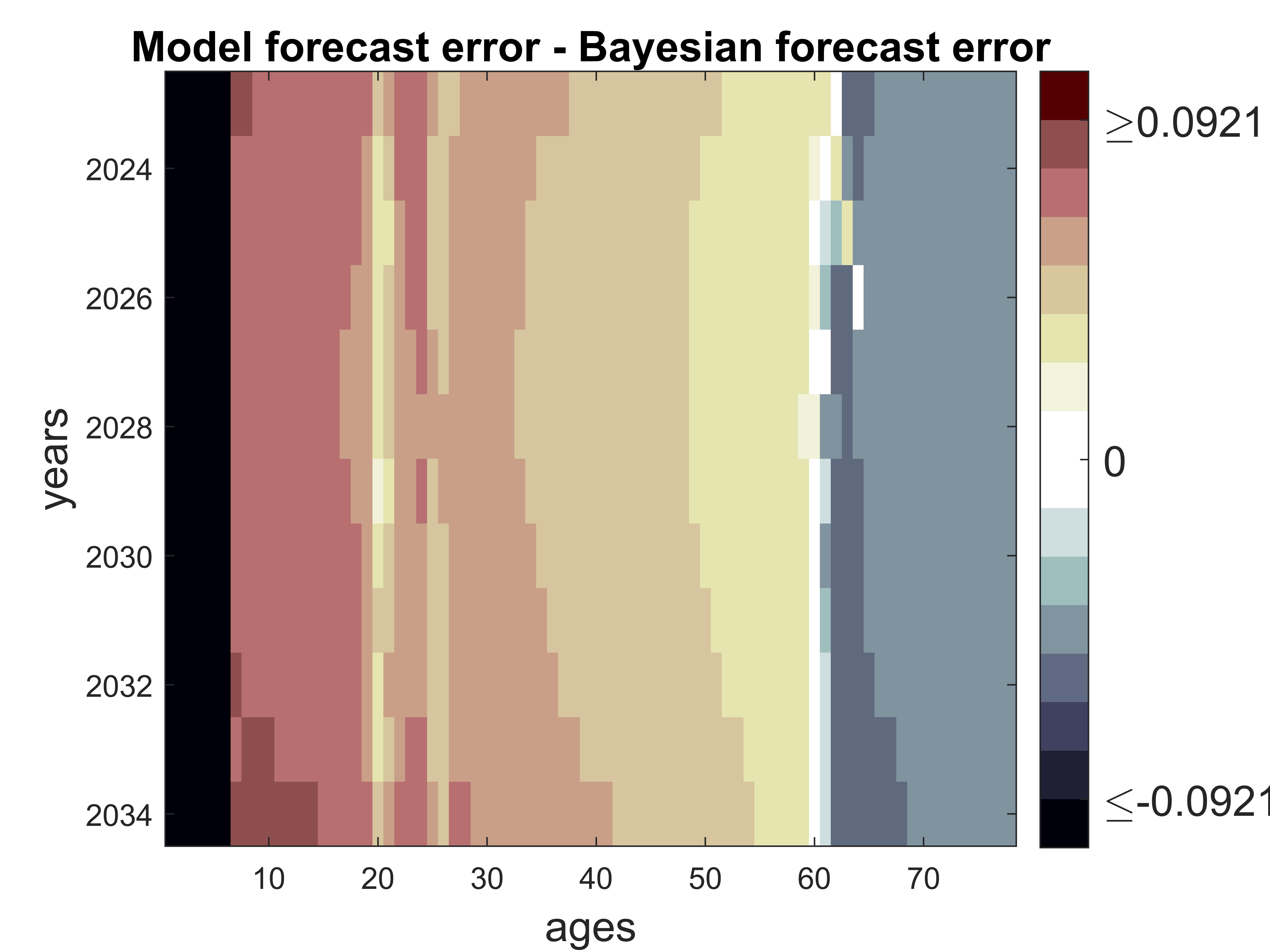}
\caption{The forecast errors between extended model and Bayesian method.}
\label{fig:forecast_error}

\end{figure}

The strength of the developed model is described in Figures \ref{fig:forecastmethods_difference} and \ref{fig:forecast_error} which gives the different in forecast results and errors between developed and Bayesian method.\\

The model consistently predicted high mortality rates for younger ages every year, with rates decreasing as ages increased. At older ages, both the model and Bayesian method produced similar mortality rates across all years. This suggests that the Bayesian method generally provided lower forecast mortality rates compared to the developed model. Furthermore, the Bayesian method demonstrated greater stability and reliability in predicting mortality rates when compared to the model.\\

Additionally, Figure~\ref{fig:forecast_error} illustrates the forecast error, showing that the developed model had very low errors for younger ages from 0--5 years, indicating its predictions are closely matched to the data. However, errors increased for ages between 10 and 60, suggesting less accuracy in this age range. 

The model became more accurate for individuals aged 60 and above, with decreased forecast errors. This pattern suggests that the developed model performed well for younger and old ages and not in the middle ages between 10--40. Furthermore, as seen in Figure~\ref{fig:forecastmethods_difference}, this reinforces that the Bayesian method consistently provided accurate and reliable predictions compared to the developed model for many years.
\section{Conclusion and Further work.}
\label{sec:conclusion}
	The increased life expectancy has brought a negative effect on the performance of several sectors, such as life insurance and pension systems. The problem comes from the use of mortality data of developed countries from lifetables because they fail to capture the decline in mortality rates across periods of time. In developing countries like Tanzania, life tables are commonly used for social and economic planning in the different sectors.\\

	In this study, we have developed a mortality rate model and projected the mortality rate using Tanzanian data. To capture the situation in limited data on mortality rates, the power exponential model that varies with both age and time-dependent was developed. The model successfully forecasted the mortality rate in Tanzania for the period of $2023-2034$.\\ 
		
   	The exhaustive numerical results showed that the developed model provided almost similar results when compared with the generated Tanzania mortality data and the existing forecasting Bayesian method, with a small forecast error between them. Moreover, the model is extremely easy to calibrate. This demonstrates that the model may be used to predict Tanzania's mortality rates.\\
   	
     In terms of usage, the model gave better results at ages above 60, which is preferable for pension systems since understanding how many people will live to a high age is necessary for pension planning. Beyond its utility for pension systems, the model can be embraced by demographers and statisticians who employ the cohort component method for population projection as one of the main input components to develop population projection results. For example, Tanzania uses the cohort component method to produce their population projections. Furthermore, as the model produces mortality rates for each year, one can use it together with actuarial concepts to generate a lifetable for many years. Additionally, the life table can be employed to generate transitional probability matrices annually, facilitating pension system projections when used with multi-state models. \\
     
     Further research directions, including the extension of the mortality rates model with many parameters to vary with both ages and time, would be an interesting point of comparison with our model. Also, a developed model with different forecasting methods, like Lee-Carter, with different data sets can be used to produce different mortality rate forecasting results.
\section*{Acknowledgement}
We are so grateful to the International Science Programme at Uppsala University for the support in the framework of the Eastern Africa Universities Mathematics Programme (EAUMP) and to the research environment in the Mathematics Departments at the University of Dar–es–Salaam and the Imperial College of London.\\
\section*{Declaration of interest}
The authors state that they have no conflict of interest related to this research work.

\bibliographystyle{tfcad}
\bibliography{Paper1_arxiv}
\end{document}